\documentclass[a4paper,reqno,12pt]{amsart}
\usepackage[a4paper,margin=1in]{geometry}
\usepackage{graphicx}
\usepackage{amsmath,amssymb,amsthm}

\theoremstyle{plain}
\newtheorem{theorem}{Theorem}[section]

\newtheorem{proposition}[theorem]{Proposition}

\theoremstyle{definition}

\numberwithin{equation}{section}

\title[On the first  bifurcation of Stokes waves]{On the first  bifurcation of Stokes waves}

\author{Vladimir Kozlov$^1$}
\address{$^1$Department of Mathematics, Link\"oping University, SE-581 83 Link\"oping, Sweden.}

\begin{document}
	
\begin{abstract} We consider Stokes water waves on the vorticity flow in a two-dimensional channel of finite depth. In the paper \cite{Koz1a} it was proved existence of subharmonic bifurcations on a branch of Stokes waves. Such bifurcations occur near the first bifurcation in the set of Stokes waves. Moreover it is shown in that paper that the bifurcating solutions build a connected continuum  containing large amplitude waves.  This fact was proved under a certain assumption concerning the second eigenvalue of the Frechet derivative. In this paper we investigate this assumption and present explicit conditions when it is satisfied.

\end{abstract}


\maketitle

\section{Formulation of the problem}

Stokes and solitary waves  were the main subject of study in the nonlinear water wave theory up to 1980.
In 1980 (see Chen \cite{Che} and Saffman \cite{Sa}) it was discovered numerically and in 2000 (see  Buffoni, Dancer and Toland, \cite{BDT1,BDT2}) this was supported theoretically for the ir-rotational case for a flow of infinite depth that there exist new types of periodic waves with several crests on the period (the Stokes wave has only one crest). These waves occur as a result of bifurcation on a branch of Stokes waves when they approach the wave of greatest amplitude.
In my papers \cite{Koz1} and \cite{Koz1a} the existence of subharmonic bifurcations was proved on branches of Stokes waves on vorticity flow.
The main result in the latest paper \cite{Koz1a} is proved under a certain assumption on the second eigenvalue of the Frechet derivative. The main goal of this paper is to study this assumption and to give an explicit conditions for its validity.

Consider steady surface waves in a two-dimensional channel bounded below by a flat,
rigid bottom and above by a free surface that does not touch the bottom. The surface tension is neglected and the water motion can be rotational.
In appropriate Cartesian coordinates $(X, Y )$, the bottom coincides with the
$X$-axis and gravity acts in the negative $Y$ -direction. We choose the frame of reference so that the velocity field is time-independent as well as the free-surface profile
which is supposed to be the graph of $Y = \xi(X)$, $x \in \Bbb R$, where $\xi$ is a positive and
continuous unknown function. Thus
$$
{\mathcal D}={\mathcal D}_\xi = \{X\in \Bbb R, 0 <Y < \xi(X)\},\;\;{\mathcal S}={\mathcal S}_\xi=\{X\in\Bbb R,\;Y=\xi(X)\}
$$
is the water domain and the free surface respectively. We will use the stream function $\Psi$, which is connected with the velocity vector $({\bf u},{\bf v})$ as ${\bf u}=c-\Psi_Y$ and ${\bf v}=\Psi_X$, where $c$ is the wave speed.

We assume that $\xi$ is a positive,  periodic function having period $\Lambda>0$ and that $\xi$   is even and strongly monotonically decreasing on the interval $(0,\Lambda/2)$.  Since the surface tension is neglected, $\Psi$ and
$\xi$ after a certain scaling satisfy the following free-boundary problem (see for example \cite{KN14}):
\begin{eqnarray}\label{K2a}
&&\Delta \Psi+\omega(\Psi)=0\;\;\mbox{in ${\mathcal D}_\xi$},\nonumber\\
&&\frac{1}{2}|\nabla\Psi|^2+\xi=R\;\;\mbox{on ${\mathcal S}_\xi$},\nonumber\\
&&\Psi=1\;\;\mbox{on ${\mathcal S}_\xi$},\nonumber\\
&&\Psi=0\;\;\mbox{for $Y=0$},
\end{eqnarray}
where $\omega\in C^{1,\alpha}$, $\alpha\in (0,1)$, is a vorticity function and $R$ is the Bernoulli constant. We assume that $\Psi$ is even, $\Lambda$-periodic in $X$ and
\begin{equation}\label{J27ba}
\Psi_Y>0\;\;\mbox{on $\overline{{\mathcal D}_\xi}$},
\end{equation}
which means that the flow is unidirectional.

The Frechet derivative for the problem is evaluated for example in \cite{KL2}, \cite{Koz1},  and the corresponding eigenvalue problem for the  Frechet derivative has the form
\begin{eqnarray}\label{J17ax}
&&\Delta w+\omega'(\Psi)w+\mu w=0 \;\;\mbox{in ${\mathcal D}_\xi$},\nonumber\\
&&\partial_\nu w-\rho w=0\;\;\mbox{on ${\mathcal S}_\xi$},\nonumber\\
&&w=0\;\;\mbox{for $Y=0$},
\end{eqnarray}
where $\nu$ is the unite outward normal to $Y=\xi(X)$ and
\begin{equation}\label{Sept17aa}
\rho=
\rho(X)=\frac{(1+\Psi_X\Psi_{XY}+\Psi_Y\Psi_{YY})}{\Psi_Y(\Psi_X^2+\Psi_Y^2)^{1/2}}\Big|_{Y=\xi(X)}.
\end{equation}
The function $w$ in (\ref{J17ax}) is supposed also to be even and $\Lambda$-periodic. 

Let us introduce several function spaces. Let $\alpha\in (0,1)$ and $k=0,1,\ldots$. The space $C^{k,\alpha}({\mathcal D})$ consists of bounded  functions in ${\mathcal D}$ such that the norms $C^{k,\alpha}({\mathcal D}_{a,a+1})$ are uniformly bounded with respect to $a\in\Bbb R$. Here
$$
{\mathcal D}_{a,a+1}=\{(X,Y)\in \overline{\mathcal D},\;:\,a\leq x\leq a+1\}.
$$
The space   $C^{k,\alpha}_{0,\Lambda }({\mathcal D})$ \big($C^{k,\alpha}_{0,\Lambda, e}({\mathcal D})\big)$ consists of $\Lambda$-periodic ($\Lambda$-periodic and even) functions, which belong to $C^{k,\alpha}({\mathcal D})$ and vanish at $Y=0$.
Similarly, we define the space  $C^{k,\alpha}_{\Lambda}(\Bbb R)$  ($C^{k,\alpha}_{\Lambda, e}(\Bbb R)$) consisting of functions in $C^{k,\alpha}(\Bbb R)$, which are $\Lambda$-periodic ($\Lambda$-periodic and even).

We will consider a branch of Stokes water waves depending on a parameter $t\geq 0$, i.e.
\begin{equation}\label{J3b}
\xi=\xi(X,t),\;\;\psi=\psi(X,Y;t),\;\;\Lambda=\Lambda(t).
\end{equation}
For each $t$ the functions $\xi\in C^{2,\alpha}_{\Lambda,e}(\Bbb R)$ and $\Psi\in C^{3,\alpha}_{\Lambda,e}({\mathcal D})$. This branch starts from a uniform stream solution for $t=0$. 
The dependence on $t$ is analytic in the sense explained in Sect. \ref{SJ29a}. The definition of uniform stream solution together with the dispersion equation which is required for existence of the branch of the Stokes waves (\ref{J3b}) is given in the next section \ref{SJ6}. Existence of such branches was a subject of many papers. In the case of non-zero vorticity we note a fundamental work Constantin and Strauss \cite{CSst}, where a bifurcation branches for the flow with vorticity was constructed for the first time. In the case with variable period we refer to the papers Kozlov and Lokharu \cite{KL1}, \cite{KL3}.

The first (lowest eigenvalue of the problem (\ref{J17ax})) is always negative and simple (see \cite{Koz1a}) and the second one we denote by $\mu(t)$.
 Assume that

\bigskip
{\bf Assumption} There exists $t_0>0$ such that $\mu(t)\geq 0$ for $t\in (0,t_0)$ and $\mu(t)<0$ for $t\in (t_0,t_0+\epsilon)$ for a certain positive $\epsilon$.

\bigskip
This assumption describes the first bifurcation point $t_0$ on the branch (\ref{J3b}) in the class of Stokes waves of period $\Lambda(t)$. It is convenient to separate two types of  bifurcations of branches of Stokes waves:

\medskip
\noindent
(i) in the class of $\Lambda(t)$-periodic solutions (Stokes bifurcation);

\noindent
(ii) in the class of $M\Lambda(t)$-periodic solutions (M-subharmonic bifurcation);

\medskip
One can consider a more general class of bifurcation in the class of bounded solution. This leads to a more challenging problem in the theory of the Hamiltonian systems with periodic coefficient. A center manifold reduction to a finite dimensional system is presented in \cite{KT}. In this paper we confine ourselves only to the cases (i) and (ii).
The following theorem is proved in \cite{Koz1a}

\begin{theorem} Let {\bf Assumption} be fulfilled. Then there exists an integer $M_0$ and pairs $(t_M,M)$, where $M$ is integer, $M>M_0$ and $t_M>t_0$, satisfying
$$
t_M\to t_0\;\;\;\mbox{as}\;\; M\to\infty,
$$
 such that $t_M$ is $M$- subharmonic bifurcation point. There are no subharmonic bifurcations for $t<t_0$.
\end{theorem}

Moreover in Theorem 9.2, \cite{Koz1a}, a structure of the set of bifurcating solutions is given. In particular it was shown that the bifurcating solutions build a connected continuum  containing large amplitude waves.

The main aim of this paper is to give explicit conditions for validity of {\bf Assumption}. Our analysis consists of two parts:

\medskip
\noindent
(i) analysis of behaviour of $\mu(t)$ for small $t$;

\noindent
(ii) analysis of $\mu(t)$ for large positive $t$.

\medskip
For $t=0$, $\Lambda(0)=\Lambda_0$ and $\mu(0)=0$. Our first goal is to study the functions $\Lambda(t)$ and $\mu(t)$ for small $t$. One of the results is the following. It's quite straightforward to show that these functions has the following asymptotic representations
$$
\mu (t)=\mu_2t^2+0(t^3)\;\;\mbox{and}\;\; \Lambda(t)=\Lambda_0+\Lambda_2t^2+O(t^3),
$$
where $\Lambda_0=\Lambda(0)$.
It is proved that
\begin{equation}\label{J3a}
\mu_2=C\Lambda_2
\end{equation}with a positive constant $C$ to be evaluated later.
To prove formula (\ref{J3a}), first we study the function
\begin{equation}\label{J7a}
\lambda(t)=\frac{\Lambda_0}{\Lambda(t)}=1+\lambda_2t^2+O(t^3)
\end{equation}
and established the relation
\begin{equation}\label{J3aa}
-4\lambda_2\tau_*^2\int_0^d\gamma(Y;\tau_*)^2dY=\mu_2\int_0^d \gamma(Y;\tau_*)^2\frac{dY}{\Psi_Y},
\end{equation}
where $\gamma(Y;\tau)$ solves the problem (\ref{Okt6b}).
Since $\Lambda_2=-\lambda_2\Lambda_0$, the last relation implies (\ref{J3a}) with a positive constant $C$. Thus the sign of $\mu_2$ is the same as of $\Lambda_2$ and opposite to the sign of $\lambda_2$.


In the  irrotational case, i.e.  $\omega=0$, we study the dependence of $\mu_2$ (actually of $\lambda_2$) on the frequency  $\tau>0$ at $t=0$ connected with the Froude number $F=d_-^{-3/2}$ by\footnote{It follows from (\ref{Ju9a})}
\begin{equation}\label{J29b}
\tau\coth\tau=\Big(\frac{F+\sqrt{F^2+8}}{4}\Big)^3F,
\end{equation}
where the right-hand side is monotone with respect to $F$.
We prove that
$$
\mu_2(\tau)>0\;\;\mbox{for all $\tau<\tau_0$, where}\;\;\tau_0\approx 1,992.
$$
In terms of the Froude number the eigenvalue $\mu(t)$ is positive when
\begin{equation}\label{Ju5b}
F<F_0,\;\;F_0\approx 1,399.
\end{equation}
This give a condition for validity of the first part in {\bf Assumption}.

Let us turn to the second part of the above assumption. It is enough to show an appearance of  negative eigenvalues of the Frechet derivative when $t\to \infty$. According to Corollary 2.2, in Kozlov and Lokharu \cite{KL1}, there exists a sequence
 $\{t_j\}$, $j=1,\ldots$, such that

 \medskip
\noindent
 a). $\xi(0,t_j)$ tends to $R$ when $j\to\infty$ (extreme wave) or

\noindent
 b). $\xi(0,t_j)$ tends to a solitary wave as $j\to\infty$

\medskip

 In the case a) the limit configuration is the extreme wave with the angle $120^{\circ}$ at the crest (see  Amick,  Fraenkel and Toland \cite{T2}, Plotnikov \cite{P2}, McLeod \cite{McL} and Varvaruca and Weiss \cite{VW1}), and the appearance of negative eigenvalues follow from Theorem 3.1, \cite{Koz1} and \cite{KL2}.




To  show that the option b) is impossible we choose parameters of the problem such that solitary waves are excluded. We will do this by using known upper estimates for the Froude number of solitary waves.

The best known upper estimate for the Froude number of solitary wave, which follows from  Starr \cite{Star} (see also  Keady and  Pritchard \cite{KP74} and Introduction of Wheeler \cite{We}), is the following
\begin{equation}\label{J29a}
F<\sqrt{2}.
\end{equation}
This means that for $F>\sqrt{2}$ there are no solitary waves. Hence every global branch of Stokes waves must approach a Stokes waves of maximal amplitude which have the angle $120^\circ$ at the crest. According to Theorem 3.1 \cite{Koz1} this fact implies appearance of infinitely many negative eigenvalues of the Frechet derivative when $t\to\infty$. 
Unfortunately the above estimate for $F$ is not enough, since $\sqrt{2}>F_0$.
Another upper estimate for the Froude number obtained numerically (see Miles \cite{Mil80}, Longuet-Higgins and  Fenton \cite{LHF74}, Hunter and  Vanden-Broeck \cite{HVB83} and Introduction in Wheeler \cite{We}) is $F<1,29$. Hence
\begin{equation}\label{J29da}
\mbox{if}\;\;1,29<F<1,399\;\;\mbox{then {\bf Assumption} is valid}.
\end{equation}
The upper estimate here is optimal but the lower estimate can be  possibly improved, since it guarantees infinitely many negative eigenvalues of the Frechet derivative when $t\to\infty$, but for the validity of the {\bf Assumption} it is sufficient to have only one negative eigenvalue.

\subsection{Uniform stream solution, dispersion equation}\label{SJ6}

The uniform stream solution $\Psi=U(Y)$ with the constant depth $\eta =d$  satisfies the problem
\begin{eqnarray}\label{X1}
&&U^{''}+\omega(U)=0\;\;\mbox{on $(0;d)$},\nonumber\\
&&U(0)=0,\;\;U(d)=1,\nonumber\\
&&\frac{1}{2}U'(d)^2+d=R.
\end{eqnarray}
In order to find solutions to this problem we introduce a parameter $s=U'(0)$. We assume that
 $s>s_0:=2\sqrt{\max_{\tau\in [0,1]}\Omega(\tau)}$, where
$$
\Omega(\tau)=\int_0^\tau \omega(p)dp.
$$
Then the problem (\ref{X1}) has a solution $(U,d)$ with a strongly monotone function $U$ for
\begin{equation}\label{M6ax}
R={\mathcal R}(s):=\frac{1}{2}s^2+d(s)-\Omega(1).
\end{equation}
The solution is given by
\begin{equation}\label{F22a}
Y=\int_0^U\frac{d\tau}{\sqrt{s^2-2\Omega(\tau)}},\;\;d=d(s)=\int_0^1\frac{d\tau}{\sqrt{s^2-2\Omega(\tau)}}.
\end{equation}
If we consider (\ref{M6ax}) as the equation with respect to $s$ then it is solvable if $R\geq R_c$, where
\begin{equation}\label{F27a}
R_c=\min_{s\geq s_0}{\mathcal R}(s),
\end{equation}
and it has two solutions if
$R\in (R_c,R_0)$, where
\begin{equation}\label{D19ba}
R_0={\mathcal R}(s_0).
\end{equation}
We denote by $s_c$ the point where the minimum in (\ref{F27a}) is attained.

Existence of small amplitude Stokes waves is determined by the dispersion equation (see, for example, \cite{KN14}). It is defined as follows.
 The strong monotonicity of $U$ guarantees that the problem
\begin{equation}\label{Okt6b}
\gamma^{''}+\omega'(U)\gamma-\tau^2\gamma=0,\;\; \gamma(0,\tau)=0,\;\;\gamma(d,\tau)=1,
\end{equation}
has a unique solution $\gamma=\gamma(y,\tau)$ for each $\tau\in\Bbb R$, which is even with respect to $\tau$ and depends analytically on $\tau$.
Introduce the function
\begin{equation}\label{Okt6ba}
\sigma(\tau)=\kappa\gamma'(d,\tau)-\kappa^{-1}+\omega(1),\;\;\kappa=U'(d).
\end{equation}
It depends also analytically on $\tau$ and it is strongly increasing with respect to $\tau>0$. Moreover it is an even function.
The dispersion  equation (see, for example \cite{KN14})  is the following
\begin{equation}\label{Okt6bb}
\sigma(\tau)=0.
\end{equation}
It has a positive solution if
\begin{equation}\label{D17a}
\sigma(0)<0.
\end{equation}
By \cite{KN14} this is equivalent to $s+d'(s)<0$ or what is the same
\begin{equation}\label{D19b}
1<\int_0^d\frac{dY}{U'^2(Y)}.
\end{equation}
The right-hand side here is equal to $1/F^2$ where $F$ is the Froude number (see \cite{We} and \cite{KLW}). Therefore (\ref{D19b}) means that $F<1$, which is well-known condition for existence of  Stokes waves of small amplitude.
Another equivalent formulation is given by requirement (see, for example \cite{KN11})
\begin{equation}\label{M3aa}
s\in (s_0,s_c).
\end{equation}
The existence of such $s$ is guaranteed by $R\in (R_c,R_0)$. One more formula for the Froude number is the following
\begin{equation}\label{Au21a}
\frac{1}{F^2(s)}=-\frac{d'(s)}{s},
\end{equation}
where the Froude number $F(s)$ corresponds to the uniform stream  solution $(U(Y;s),d(s))$ and $R={\mathcal R}(s)$.
Therefore
$$
{\mathcal R}'(s)=s(1-F^{-2}(s)).
$$
The value $\sigma(0)$ admits the following representation (see \cite{KN14}):
$$
\sigma(0)=-\frac{3}{2\kappa}\frac{{\mathcal R}'(s)}{d'(s)}=\frac{3(F^2(s)-1)}{2\kappa},
$$
where we have used (\ref{Au21a}) to verify the second equality above.

The function $\sigma$ has the following asymptotic representation
$$
\sigma(\tau)=\kappa\tau +O(1)\;\;\mbox{for large $\tau$}
$$
and equation (\ref{Okt6bb}) has a unique positive root, which will be denoted by $\tau_*$. It is connected with $\Lambda_0$ by the relation
\begin{equation}\label{Ju8a}
\tau_*=\frac{2\pi}{\Lambda_0}.
\end{equation}

To give another representation of the function $\sigma$ we introduce
\begin{equation}\label{F28a}
\rho_0=\frac{1+U'(d)U^{''}(d)}{U'(d)^2}
\end{equation}
and note that
$$
\frac{1+U'(d)U^{''}(d)}{U'(d)^2}=\kappa^{-2}-\frac{\omega(1)}{\kappa}.
$$
Hence another form for (\ref{Okt6ba}) is
\begin{equation}\label{M21aa}
\sigma(\tau)=\kappa\gamma'(d,\tau)-\kappa\rho_0.
\end{equation}

The following problem will be used in the asymptotic analysis of the branch (\ref{J3b}) for small $t$:
\begin{eqnarray}\label{J7b}
&&v^{''}+\omega'(U)v-\tau^2v=f\;\;\mbox{on $(0,d)$},\nonumber\\
&&v'(d)-\rho_0v(d)=g\;\;\mbox{and}\;\;v(0)=0.
\end{eqnarray}
\begin{proposition} Let $\tau\geq 0$ and $\tau\neq\tau_*$. Let also $f\in C^{1,\alpha}([0,d])$ and $g$ be a constant. Then the problem
{\rm (\ref{J7b})} has a unique solution $v\in C^{3,\alpha}$.
If $\tau=\tau_*$ then the problem {\rm (\ref{J7b})} has  one dimensional kernel which consists of function
$$
c\gamma(Y;\tau_*).
$$
\end{proposition}

\section{A connection between the functions $\mu(t)$ and $\Lambda(t)$ for small $t$}

In this section we prove formula (\ref{J3aa}). It appears that the partial hodograph transform is very useful for this purpose.

\subsection{Partial hodograph transform}\label{SJ29a}

In what follows we will study  branches of Stokes waves $(\Psi(X,Y;t),\xi(X;t))$ of period $\Lambda(t)$, $t\geq 0$, started from the uniform stream  at $t=0$. The existence of such branches is established  first in \cite{CSst} with fixed period but variable $R$ and in \cite{KL1} for variable $\Lambda$ and fixed $R$. In our case of variable $\Lambda$
it is convenient to make the following change of variables
\begin{equation}\label{D11a}
x=\lambda X,\;\;y=Y,\;\;\lambda=\frac{\Lambda_0}{\Lambda(t)}
\end{equation}
in order to deal with the problem with a fixed period. Here as before
$
\Lambda_0:=\Lambda(0)=2\pi/\tau_*,
$
where $\tau_*$ is the root of the equation  (\ref{Okt6bb}).
As the result we get
\begin{eqnarray}\label{Okt6aa}
&&\Big(\lambda^2\partial_x^2+\partial_y^2\Big)\psi+\omega(\psi)=0\;\;\mbox{in $D_\eta$},\nonumber\\
&&\frac{1}{2}\Big(\lambda^2\psi_x^2+\psi_y^2\Big)+\eta=R\;\;\mbox{on $B_\eta$},\nonumber\\
&&\psi=1\;\;\mbox{on $B_\eta$},\nonumber\\
&&\psi=0\;\;\mbox{for $y=0$},
\end{eqnarray}
where
$$
\psi(x,y;t)=\Psi(\lambda^{-1}x,y;t)\;\;\mbox{and}\;\;\eta(x;t)=\xi(\lambda^{-1} x;t).
$$
Here all functions have the  same period $\Lambda_0$, $D_\eta$ and $B_\eta$  are the domain and the free surface  after the change of variables (\ref{D11a}).

From (\ref{J27ba}) it follows that
$$
\psi_y>0\;\;\mbox{in $\overline{D_\eta}$}.
$$
Using the change of variables
$$
q=x,\;\;p=\psi,
$$
we get
$$
q_x=1,\;\;q_y=0,\;\;p_x=\psi_x,\;\;p_y=\psi_y,
$$
and
\begin{equation}\label{F28b}
\psi_x=-\frac{h_q}{h_p},\;\;\psi_y=\frac{1}{h_p},\;\;dxdy=h_pdqdp.
\end{equation}

System (\ref{Okt6aa}) in the new variables takes the form
\begin{eqnarray}\label{J4a}
&&\Big(\frac{1+\lambda^2h_q^2}{2h_p^2}+\Omega(p)\Big)_p-\lambda^2\Big(\frac{h_q}{h_p}\Big)_q=0\;\;\mbox{in $Q$},\nonumber\\
&&\frac{1+\lambda^2h_q^2}{2h_p^2}+h=R\;\;\mbox{for $p=1$},\nonumber\\
&&h=0\;\;\mbox{for $p=0$}.
\end{eqnarray}
Here
$$
Q=\{(q,p)\,:\,q\in\Bbb R\,,\;\;p\in (0,1)\}.
$$
The uniform stream solution corresponding to the solution $U$ of (\ref{X1}) is
\begin{equation}\label{M4c}
H(p)=\int_0^p\frac{d\tau}{\sqrt{s^2-2\Omega(\tau)}},\;\;s=U'(0)=H_p^{-1}(0).
\end{equation}
One can check that
\begin{equation}\label{J18aaa}
H_{pp}-H_p^3\omega(p)=0
\end{equation}
or equivalently
\begin{equation}\label{J18aa}
\Big(\frac{1}{2H_p^2}\Big)_p+\omega(p)=0.
\end{equation}
Moreover it satisfies the boundary conditions
\begin{equation}\label{M4ca}
\frac{1}{2H_p^2(1)}+H(1)=R,\;\;H(0)=0.
\end{equation}
The Froude number in new variables can be written as
$$
\frac{1}{F^2}=\int_0^1H_p^3dp.
$$

 Then according to Theorem 2.1, \cite{KL1} there exists a branch of solutions to (\ref{J4a})
\begin{equation}\label{J4ac}
h=h(q,p;t):[0,\infty)\rightarrow C^{2,\gamma}_{pe}(\overline{Q}),\;\;\lambda=\lambda(t):[0,\infty)\rightarrow (0,\infty),
\end{equation}
which  has a real analytic reparametrization locally around each $t\geq 0$.

\subsection{Bifurcation equation}

In order to find bifurcation points and bifuracating solutions we put $h+w$ instead of $h$ in (\ref{J4a}) and introduce the operators
\begin{eqnarray*}
&&{\mathcal F}(w;t)=\Big(\frac{1+\lambda^2(h_q+w_q)^2}{2(h_p+w_p)^2}\Big)_p
-\Big(\frac{1+\lambda^2h_q^2}{2h_p^2}\Big)_p\\
&&-\lambda^2\Big(\frac{h_q+w_q}{h_p+w_p}\Big)_q+\lambda^2\Big(\frac{h_q}{h_p}\Big)_q
\end{eqnarray*}
and
$$
{\mathcal G}(w;t)=\frac{1+\lambda^2(h_q+w_q)^2}{2(h_p+w_p)^2}-\frac{1+\lambda^2h_q^2}{2h_p^2}+w
$$
acting on $\Lambda_0$-periodic, even functions $w$ defined in $Q$. After some cancelations we get
$$
{\mathcal F}={\mathcal J}_p+{\mathcal I}_q,\;\;{\mathcal G}={\mathcal J}+w,
$$
where
$$
{\mathcal J}={\mathcal J}(w;t)=\frac{\lambda^2h_p^2(2h_q+w_q)w_q-(2h_p+w_p)(1+\lambda^2h_q^2)w_p}{2h_p^2(h_p+w_p)^2}
$$
and
$$
{\mathcal I}={\mathcal I}(w;t)=-\lambda^2\frac{h_pw_q-h_qw_p}{h_p(h_p+w_p)}.
$$
Both these functions are well defined  for small $w_p$.
Then the problem for finding solutions close to $h$ is the following
\begin{eqnarray}\label{F19a}
&&{\mathcal F}(w;t)=0\;\;\mbox{in $Q$}\nonumber\\
&&{\mathcal G}(w;t)=0\;\;\mbox{for $p=1$}\nonumber\\
&&w=0\;\;\mbox{for $p=0$}.
\end{eqnarray}

Furthermore, the Frechet derivative (the linear approximation of the functions ${\mathcal F}$ and ${\mathcal G}$) is the following
\begin{equation}\label{J4aa}
Aw=A(t)w=\Big(\frac{\lambda^2h_qw_q}{h_p^2}-\frac{(1+\lambda^2h_q^2)w_p}{h_p^3}\Big)_p-\lambda^2\Big(\frac{w_q}{h_p}-\frac{h_qw_p}{h_p^2}\Big)_q
\end{equation}
and
\begin{equation}\label{J4aba}
{\mathcal N}w={\mathcal N}(t)w=(N w-w)|_{p=1},
\end{equation}
where
\begin{equation}\label{J4ab}
N w=N(t)w=\Big(-\frac{\lambda^2h_qw_q}{h_p^2}+\frac{(1+\lambda^2h_q^2)w_p}{h_p^3}\Big)\Big|_{p=1}.
\end{equation}
The eigenvalue problem for the Frechet derivative, which is important for the analysis of bifurcations of the problem
(\ref{F19a}), is the following
\begin{eqnarray}\label{M1a}
&&A(t)w=\mu w\;\;\mbox{in $Q$},\nonumber\\
&&{\mathcal N}(t)w=0\;\;\mbox{for $p=1$},\nonumber\\
&&w=0\;\;\mbox{for $p=0$}.
\end{eqnarray}

For $t=0$ and $\mu=0$ this problem becomes
\begin{eqnarray}\label{J15a}
&&A_0w:=-\Big(\frac{w_{p}}{H_p^3}\Big)_p-\Big(\frac{w_{q}}{H_p}\Big)_q=0\;\;\mbox{in $Q$},\nonumber\\
&&B_0w:=-\frac{w_{p}}{H_p^3}+w=0\;\;\mbox{for $p=1$},\nonumber\\
&&w=0\;\;\mbox{for $p=0$}.
\end{eqnarray}
Since the function $H$ depends only on $p$ this problem admits the separation of variables and its solutions are among the functions
\begin{equation}\label{J15aa}
v(q,p)=\alpha(p)\cos (\tau q),\;\;\mbox{where}\;\;\tau=k\tau_*,\;\;k=0,1,\ldots.
\end{equation}
According to \cite{Koz1} the function (\ref{J15aa}) solves (\ref{J15a}) if and only if
$$
\alpha(p)=\gamma(H(p);\tau)H_p,
$$
where the function $\gamma(Y;\tau)$ solves the euation (\ref{Okt6b}) and $\sigma(\tau)=0$. Therefore if $\tau\neq \tau_*$ then the problem (\ref{J15a}) has no non-trivial solutions. If $\tau=\tau_*$ then
the kernel of the above operator is one dimensional in the class of $\Lambda_0:=2\pi/\tau_*$ periodic, even function and it is given by
\begin{equation*}
v=\alpha(p)\cos(\tau_*q),\;\;  \alpha(p)=\gamma(H(p);\tau_*)H_p.
\end{equation*}


We will need also the problem
\begin{eqnarray}\label{M6a}
&& -\Big(\frac{u_{p}}{H_p^3}\Big)_p+\frac{\tau^2u}{H_p}=F\;\;\mbox{on (0,d)}\nonumber\\
&&u(0)=0,\;\;-\frac{u_{p}}{H_p^3}+u=c\;\;\mbox{for $p=1$},
\end{eqnarray}
where $F\in C^{0,\alpha}([0,1])$ and $c$ is a constant.
Clearly this problem is
elliptic and uniquely solvable for all $\tau\geq 0$, $\tau\neq\tau_*$, the problem (\ref{M6a}) has a unique solution in  $C^{2,\alpha}([0,1])$. This solution is given by
$$
u(p)=v(H(p))H_p(p),
$$
where $v(Y)$ solves the problem (\ref{J7b}) with $f=F(H(y))$ and $g=c$.

\subsection{Stokes waves for small $t$}
Here we consider asymptotics of solutions of (\ref{J4ac}) for small $t$. For this purpose we take
$$
h=H(p)
$$
and represent the solution in the form
$$
H(p)+w(q,p,t),\;\;w=tv,
$$
where 
\begin{equation}\label{Okt11b}
v(q,p;t)=v_0(q,p)+tv_1(q,p)+t^2v_2(q,p)+\cdots
\end{equation}
The function $\lambda=\lambda(t)$ is sought in the form
\begin{equation}\label{M5a}
\lambda(t)=1+\lambda_2t^2+O(t^4).
\end{equation}
The coefficients $\lambda_1$ and $\lambda_3$ in the above formula are zero as one can easily see from the forthcoming calculations.
Our aim is to find Stokes waves close to $H$. Since the functions $w$, $v$ and $\lambda$ analytically depend on $t$ it is sufficient to find coefficients $v_j$ and $\lambda_j$.

In this case
$$
{\mathcal J}=A_1\Big(1+\frac{w_p}{H_p}\Big)^{-2}+A_2\Big(1+\frac{w_p}{H_p}\Big)^{-2},
$$
where
$$
A_1=-\frac{w_p}{H_p^3}
$$
and
$$
A_2=\frac{\lambda^2w_q^2}{2H_p^2}-\frac{w_p^2}{2H_p^4}.
$$
Therefore
$$
{\mathcal J}={\mathcal J}_1+{\mathcal J}_2+{\mathcal J}_3+O(t^4),
$$
where
$$
{\mathcal J}_1=A_1,
$$
$$
{\mathcal J}_2=A_2-2\frac{w_p}{H_p}A_1=\frac{\lambda^2w_q^2}{2H_p^2}+\frac{3}{2}\frac{w_p^2}{H_p^4}
$$
and
$$
{\mathcal J}_3=3\frac{w_p^2}{H_p^2}A_1-2\frac{w_p}{H_p}A_2=-2\frac{w_p^3}{H_p^5}-\frac{w_pw_q^2}{H_p^3}.
$$
Furthermore
$$
{\mathcal I}=-\lambda^2\frac{w_q}{H_p}\Big(1+\frac{w_p}{H_p}\Big)^{-1}
={\mathcal I}_1+{\mathcal I}_2+{\mathcal I}_3+O(t^4).
$$
Here
$$
{\mathcal I}_1=-\lambda^2\frac{w_q}{H_p},\;\;{\mathcal I}_2=\lambda^2\frac{w_qw_p}{H^2_p},\;\;{\mathcal I}_3(w)=-\lambda^2\frac{w_qw_p^2}{H^3_p}.
$$


Inserting (\ref{Okt11b}) and (\ref{M5a}) into (\ref{F19a}) and
equating terms of the same power with respect to $t$, we get
\begin{eqnarray*}
&&A_0v_0=-\Big(\frac{v_{0p}}{H_p^3}\Big)_p-\Big(\frac{v_{0q}}{H_p}\Big)_q=0\;\;\mbox{in $Q$},\\
&&B_0v_0=-\frac{v_{0p}}{H_p^3}+v_0=0\;\;\mbox{for $p=1$},\\
&&v_0=0\;\;\mbox{for $p=0$}.
\end{eqnarray*}
As we have shown in previous section the kernel of the above operator is one dimensional  and is generated by the function
\begin{equation}\label{Okt12a}
v_0=\alpha_0(p)\cos(\tau_*q),\;\;\alpha_0=\gamma(H(p);\tau_*)H_p.
\end{equation}

The next term in the asymptotics satisfies the boundary value problem
\begin{eqnarray}\label{Okt13a}
&&A_0v_1+\Big(\frac{v_{0q}^2}{2H_p^2}+\frac{3}{2}\frac{v_{0p}^2}{H_p^4}\Big)_p+\Big(\frac{v_{0q}v_{0p}}{H^2_p}\Big)_q=0\;\;\mbox{in $Q$},\nonumber\\
&&B_0v_{1}+\frac{v_{0q}^2}{2H_p^2}+\frac{3}{2}\frac{v_{0p}^2}{H_p^4}=0 \;\;\mbox{for $p=1$},\nonumber\\
&&v_1=0\;\;\mbox{for $p=0$}.
\end{eqnarray}
The solution of this problem, orthogonal to $v_0$ in $L^2$, is given by
\begin{equation}\label{Okt12aa}
v_1=\alpha_1(p)+\beta_1(p)\cos(2\tau_* q),
\end{equation}
where $\alpha_1$ and $\beta_1$ satisfy the problem (\ref{M6a}) with $\tau=0$ and $\tau=2\tau_*$ respectively with certain right-hand sides.
 Further, the term $v_2$ is found from the following problem
\begin{eqnarray*}
&&\!\!A_0v_2\!+\!\Big(\frac{v_{0q}v_{1q}}{H_p^2}\!+\!\frac{3v_{0p}v_{1p}}{H_p^4}\!+\!{\mathcal J}_3(v_0)\Big)_p\!
+\!\Big(\frac{v_{1q}v_{0p}+v_{0q}v_{1p}}{H^2_p}\!+\!{\mathcal I}_3(v_0)\Big)_q\!=\!2\lambda_2\Big(\frac{v_{0q}}{H_p}\Big)_q\;\;\mbox{on $Q$},\\
&&B_0v_2+\frac{v_{0q}v_{1q}}{H_p^2}+\frac{3v_{0p}v_{1p}}{H_p^4}+{\mathcal J}_3(v_0)=0\;\;\mbox{for $p=1$}\\
&&v_2(q,0)=0.
\end{eqnarray*}

The solvability condition for the last problem has the form
\begin{eqnarray}\label{Okt4a}
&&2\lambda_2\int_\Omega\frac{v_{0q}^2}{H_p}dqdp-\int_\Omega \Big(\Big(\frac{v_{0q}v_{1q}}{H_p^2}+\frac{3v_{0p}v_{1q}}{H_p^4}\Big)v_{0p}+\frac{v_{0q}v_{1p}+v_{1q}v_{0p}}{H^2_p}v_{0q}\Big)dqdp\nonumber\\
&&+\int_\Omega\Big(\Big(\frac{2v_{0p}^3}{H_p^5}+\frac{v_{0p}v_{0q}^2}{H_p^3}\Big)v_{0p}+\frac{v_{0p}^2v_{0q}}{H_p^3}v_{0q}\Big)dqdp=0,
\end{eqnarray}
where
$$
\Omega=\{(q,p)\in Q\,:\,-\Lambda_0/2<q<\Lambda_0/2\}.
$$
Relation (\ref{Okt4a}) can be used to find $\lambda_2$. It is quite difficult to find the sign of $\lambda_2$ from this relation but it implies a continuity of $\lambda_2$ with respect to $R$ and $\omega$. The function $v_2$ has the form
\begin{equation}\label{J12a}
v_2=\alpha_2(p)\cos(\tau_* q)+\beta_2(p)\cos(3\tau_* q),
\end{equation}
where $\alpha_2$ and $\beta_2$ satisfy the problem (\ref{M6a}) with $\tau=\tau_*$ and $\tau=3\tau_*$ respectively with certain right-hand sides.

Thus we have shown that $\lambda$ and $v$ have the form (\ref{M5a}) and (\ref{Okt11b}) respectively. More exactly $v_0$ is given by (\ref{Okt12a}), $v_1$ is represented as (\ref{Okt12aa}) and $v_2$ by (\ref{J12a}).

\subsection{Formula for $\lambda_2$ and the proof of the relation (\ref{J3aa})}

Using the representations (\ref{J4aa}) and (\ref{J4aba}) with $h=H+w$, where $w$ is evaluated in the previous section, we can write the Frechet derivative of the operators ${\mathcal J} U$ and ${\mathcal I} U$ in the form

$$
d{\mathcal J}(U)=-\frac{U_p}{H_p^3}+\Big(\frac{w_qU_q}{H_p^2}+3\frac{w_pU_p}{H_p^4}\Big)-6\frac{w_p^2U_p}{H_p^5}
-\frac{w_q^2U_p+2w_pw_qU_q}{H_p^3}+O(t^3)
$$
and
$$
d{\mathcal I}(U)=-\lambda^2\frac{U_q}{H_p}+\frac{w_pU_q+w_qU_p}{H_p^2}-\frac{w_p^2U_q+2w_qw_pU_p}{H_p^3}+O(t^3).
$$
The eigenvalue problem is described by the boundary value problem
\begin{eqnarray*}
&&(d{\mathcal J}(U))_p+(d{\mathcal I}(U))_q=(\mu_2t^2+O(t^3))U\;\;\mbox{in $Q$}\\
&&d{\mathcal J}(U)+U=0\;\;\mbox{for $p=1$}\\
&&U=0\;\;\mbox{for $p=0$.}
\end{eqnarray*}
We are looking for the eigenfunction $U$ in the form
$$
U=U(q,p;t)=U_0(q,p)+tU_1(q,p)+t^2U_2(q,p)+O(t^3),\;\;U_0=v_0.
$$
Equating terms of the same order with respect to $t$, we get
\begin{eqnarray*}
&&A_0U_1+\Big(\frac{v_{0q}U_{0q}}{H_p^2}+3\frac{v_{0p}U_{0p}}{H_p^4}\Big)_p+\Big(\frac{w_{0p}U_{0q}+v_{0q}U_{0p}}{H_p^2}\Big)_q=0\;\;\mbox{in $Q$},\\
&&B_0U_1+\Big(\frac{v_{0q}U_{0q}}{H_p^2}+3\frac{v_{0p}U_{0p}}{H_p^4}\Big)=0\;\;\mbox{for $p=1$},\\
&&U_1=0\;\;\mbox{for $p=0$}.
\end{eqnarray*}
Comparing this problem with (\ref{Okt13a}) and using that $U_0=v_0$, we conclude that $U_1=2v_1$.

Next, we write the equation for $U_2$
\begin{eqnarray}\label{Okt14a}
&&-2\lambda_2\Big(\frac{U_{0q}}{H_p}\Big)_q+A_0U_2+\Big(\frac{v_{1q}U_{0q}+v_{0q}U_{1q}}{H_p^2}+3\frac{v_{1p}U_{0p}+v_{0p}U_{1p}}{H_p^4}\Big)_p\nonumber\\
&&+\Big(\frac{v_{1p}U_{0q}+v_{1q}U_{0p}+v_{0p}U_{1q}+v_{0q}U_{1p}}{H_p^2}\Big)_q
-\Big(6\frac{v_{0p}^2U_{0p}}{H_p^5}\nonumber\\
&&+\frac{v_{0q}^2U_{0p}+2v_{0p}v_{0q}U_{0q}}{H_p^3}\Big)_p-\Big(\frac{v_{0p}^2U_{0q}+2v_{0q}v_{0p}U_{0p}}{H_p^3}\Big)_q=\mu_2U_0\;\;\mbox{in $Q$}
\end{eqnarray}\label{Okt14aa}
and the boundary equations $U_2=0$ for $p=0$ and
\begin{eqnarray*}
&&B_0U_2+\Big(\frac{v_{1q}U_{0q}+v_{0q}U_{1q}}{H_p^2}+3\frac{v_{1p}U_{0p}+v_{0p}U_{1p}}{H_p^4}\Big)\nonumber\\
&&-\Big(6\frac{v_{0p}^2U_{0p}}{H_p^5}+\frac{v_{0q}^2U_{0p}+2v_{0p}v_{0q}U_{0q}}{H_p^3}\Big)=0\;\;\mbox{for $p=1$}
\end{eqnarray*}
Since $U_0=v_0$ and $U_1=2v_1$, the solvability condition for (\ref{Okt14a}) has the form
\begin{eqnarray}\label{KK1}
&&2\lambda_2\int_{\Omega}\frac{v_{0q}^2}{H_p}dqdp-3\int_{Q_p}\Big(\frac{v_{1q}v_{0q}}{H_p^2}+3\frac{v_{0p}v_{1q}}{H_p^4}\Big)v_{0p}dqdp\\
&&-3\int_{\Omega}\frac{v_{1p}v_{0q}+v_{1q}v_{0p}}{H_p^2}v_{0q}dqdp+3\Big(\int_{Q_p} \Big(2\frac{v_{0p}^3}{H_p^5}+\frac{v_{0q}^2v_{0p}}{H_p^3}\Big)v_{0p}+\frac{v_{0p}^2v_{0q}^2}{H_p^3}\Big)dqdp\nonumber\\
&&=\mu_2\int_{\Omega} v_{0}^2dqdp.
\end{eqnarray}
Taking the sum of (\ref{KK1}) and (\ref{Okt4a}) with the factor $-3$, we get
$$
-4\lambda_2\int_\Omega\frac{v_{0q}^2}{H_p}dqdp=\mu_2\int_\Omega v_{0}^2dqdp,
$$
which coincides with (\ref{J3aa}).

\section{The coefficient $\lambda_2$ for the irrotational flow}

In this section we evaluate the coefficient $\lambda_2$
in the case $\omega=0$. The problem  (\ref{X1}) is solvable if $R\geq R_c$, where $R_c=3/2$. If $R>R_c$ then the equation
$$
\frac{1}{d^2}+2d=2R
$$
has exactly two solutions $0<d_-<1<d_+$ which are called  supercritical and subcritical, respectively. The Stokes branches appears only for the stream solutions- $(Y/d_+,d_+)$.

We will make the following change of variables
$$
X=\frac{x}{d_+},\;\;Y=\frac{y}{d_+}-1,\;\;\xi(X)=\frac{\eta(x)}{d_+}-1,\;\;\Psi(X,Y)=\psi(x,y).
$$
Then the problem (\ref{K2a}) takes the form
\begin{eqnarray}\label{Okt17a}
&&\Delta_{x,y} \psi=0\;\;\mbox{in $D_\eta$},\nonumber\\
&&|\nabla_{x,y}\psi|^2+2\theta\eta=1\;\;\mbox{on ${\mathcal S}_\eta$},\nonumber\\
&&\psi=1\;\;\mbox{for $y=\eta(x)$},\nonumber\\
&&\psi=0\;\;\mbox{for $y=-1$},
\end{eqnarray}
where
$$
\theta=d_+^3\;\;\mbox{and}\;\;D=D_\eta=\{(x,y)\,:\,x\in\Bbb R,\;-1<y<\eta(x)\}.
$$
So $\theta$ is the only parameter in the problem and $\theta\in (1,\infty)$.

We note the function $\lambda=\lambda(t)$ will be the same after the above change of variables and the same is true for  the coefficient $\lambda_2$. In the irrotational case one can derive an explicit equation for the coefficient $\lambda_2$, what will be done  below. 

\subsection{Asymptotic analysis}

In order to deal with the problem with constant period $\Lambda_0$, we perform one more change of variable
$$
s=\lambda x,\;\;\lambda=\frac{\Lambda_0}{\Lambda}.
$$
Then the problem (\ref{Okt17a}) becomes
\begin{eqnarray}\label{Okt31a}
&&(\lambda^2\partial_s^2+\partial_y^2)\psi=0\;\;\mbox{in $D$},\nonumber\\
&&\psi(s,-1)=0,\;\;\psi(s,\eta(s))=1\;\;\mbox{for $s\in\Bbb R$}\nonumber\\
&&|\partial_y\psi|^2+\lambda^2|\partial_s\psi|^2+2\theta\eta=1\;\;\mbox{for $y=\eta(s)$}.
\end{eqnarray}
We are looking for the functions $\eta$, $\psi$ and $\lambda$ in the form
\begin{equation}\label{Okt31b}
\eta=\eta(s;t))=t\eta_1(s)+t^2\eta_2(s)+t^3\eta_3(s)+\cdots,
\end{equation}
where
$$
\eta_1(s)=\cos(\tau_* s),\;\;\eta_2(s)=\beta_1+\beta_2\cos(2\tau_*s),
$$

\begin{equation}\label{Okt31ba}
\psi=\psi(s,y;t)=1+y+t\psi_1(s,y)+t^2\psi_2(s,y)+t^3\psi_3(s,y)+\cdots,
\end{equation}
$$
\psi_1(s,y)= a\cos(\tau_* s)\sinh(\tau_*(y+1)),\;\;\psi_2(s,y)=\alpha_1(y+1)+\alpha_2\cos(2\tau_* s)\sinh(2\tau_*(y+1))),
$$
and
\begin{equation}\label{Okt31bc}
\lambda=\lambda(t)=1+\lambda_2t^2+\cdots.
\end{equation}
Here $\tau_*$ is the solution of the dispersion equation
\begin{equation}\label{Okt31ca}
\nu(\tau)=\theta,\;\;\mbox{where $\nu(\tau)=\tau\coth(\tau)$,}
\end{equation}
and
\begin{equation}\label{Okt31cb}
a=-\frac{1}{\sinh\tau_*}.
\end{equation}
This implies, in particular that $\psi_1(s,0)=-\eta_1(s)$.

Inserting (\ref{Okt31b})--(\ref{Okt31bc}) into (\ref{Okt31a}), we get
\begin{equation}\label{Nov1a}
(\partial_s^2+\partial_y^2)\psi_3+2\lambda_2\partial_s^2\psi_1=0\;\;\mbox{for $s\in\Bbb R$ and $y\in (-1,0)$},
\end{equation}
\begin{equation}\label{Nov3ba}
\psi_3(s,-1)=0,
\end{equation}
\begin{equation}\label{Nov3bb}
\eta+t\psi_1(s,\eta)+t^2\psi_2(s,\eta)+t^3\psi_3(s,0)=O(t^4)
\end{equation}
and
\begin{eqnarray}\label{Nov3bc}
&&2t(\psi_{1y}(s,\eta)+t\psi_{2y}(s,\eta)+t^2\psi_{3y}(s,0))+t^2(\psi_{1y}(s,\eta)+t\psi_{2y}(s,\eta))^2\nonumber\\
&&+t^2(\psi_{1s}(s,\eta)+t\psi_{2s}(s,\eta))^2+2\theta\eta=O(t^4)
\end{eqnarray}
Direct calculations give the above expressions for $\eta_1$ and $\psi_1$ together with formulas (\ref{Okt31ca}) and  (\ref{Okt31cb}).

Let us find terms $\eta_2$ and $\psi_2$.  Using (\ref{Nov3bb}), we get the equation
$$
\beta_1+\beta_2\cos(2\tau_*s)+\alpha_1+\alpha_2\cos(2\tau_* s)\sinh(2\tau_*)+a\tau_*\cos^2(\tau_* s)\cosh(\tau_*)=0,
$$
which leads to
$$
\beta_1+\alpha_1+\frac{a}{2}\tau_*\cosh(\tau_*)=0,\;\;\beta_2+\alpha_2\sinh(2\tau_*)+\frac{a}{2}\tau_*\cosh(\tau_*)=0
$$
or due to the dispersion equation
$$
\beta_1+\alpha_1-\frac{\theta}{2}=0,\;\;\beta_2+\alpha_2\sinh(2\tau_*)-\frac{\theta}{2}=0.
$$

 From the Bernoulli relation (\ref{Nov3bc}) we derive
 $$
 2(\psi_{1yy}(s,0)\eta_1+\psi_{2y}(s,0))+\psi_{1y}^2(s,0)+\psi_{1s}^2(s.0)+2\theta\eta_2=0.
 $$
 or
 \begin{eqnarray*}
 &&2(a\tau_*^2\cos^2\tau_*s\sinh\tau_*+\alpha_1+2\alpha_2\tau_*\cos2\tau_* s\cosh(2\tau_*))+a^2\tau_*^2\cos^2\tau_*s\cosh^2\tau_*\\
 &&+a^2\tau_*^2\sin^2(\tau_* s)\sinh^2\tau_*+2\theta(\beta_1+\beta_2\cos(2\tau_*s))=0.
  \end{eqnarray*}
  Using the formula for $a$ and the dispersion equation, we get
 \begin{eqnarray*}
 &&2(-\tau_*^2\cos^2\tau_*s+\alpha_1+2\alpha_2\tau_*\cos2\tau_* s\cosh(2\tau_*))+\theta^2\cos^2\tau_*s\\
 &&+\tau_*^2\sin^2(\tau_* s)+2\theta(\beta_1+\beta_2\cos(2\tau_*s))=0.
  \end{eqnarray*}
This implies
$$
2\alpha_1+2\theta\beta_1+\frac{1}{2}(\theta^2-\tau_*^2)=0
$$
and
$$
4\tau_*\alpha_2\cosh 2\tau_*+2\theta\beta_2+\frac{1}{2}(\theta^2-3\tau_*^2)=0.
$$
Thus
$$
2(\theta-1)\beta_1+\theta+\frac{1}{2}(\theta^2-\tau_*^2)=0
$$
and
$$
(2\theta-2\nu(2\tau_*))\beta_2+\theta\nu(2\tau_*)+\frac{1}{2}(\theta^2-3\tau_*^2)=0.
$$


Let us turn to the terms of order $O(t^3)$. The equation (\ref{Nov1a}) becomes
\begin{equation}\label{Okt31caz}
(\partial_s^2+\partial_y^2)\psi_3+2\lambda_2\partial_y^2\psi_1=0,
\end{equation}
Furthermore,
\begin{equation}\label{Okt31caa}
\psi_3(s,-1)=0
\end{equation}
 and the remaining boundary relations (\ref{Nov3bb}) and (\ref{Nov3bc}) take the form
$$
\eta_3+\psi_3(s,0)+F_1=0,\;\;F_1=\psi_{2y}(s,0)\eta_1+\psi_{1y}(s.0)\eta_2+\frac{1}{2}\psi_{1yy}(s,0)\eta_1^2
$$
and
$$
2\psi_{3y}+2\theta\eta_3+F_2=0
$$
with
\begin{eqnarray*}
&&F_2=2\psi_{2yy}(s,0)\eta_1+\psi_{1yyy}(s,0)\eta_1^2+2\psi_{1yy}(s,0)\eta_2\\
&&+2\psi_{1y}(s,0)(\psi_{2y}(s,0)+\psi_{1yy}(s,0)\eta_1)+2\psi_{1s}(s,0)(\psi_{1sy}(s,0)\eta_1+\psi_{2s}(s,0)).
\end{eqnarray*}
These two equalities imply
\begin{equation}\label{Okt31cba}
2\psi_{3y}-2\theta\psi_3-2\theta F_1+F_2=0.
\end{equation}
We consider the problem (\ref{Okt31caz}), (\ref{Okt31caa}), (\ref{Okt31cba}) as the problem with respect to $\psi_3$. It has a kernel consisting of the function $\psi_1$. Therefore the solvability condition has the form
\begin{equation*}
\int_{-\Lambda_0/2}^{\Lambda_0/2}(2\theta F_1-F_2)\psi_1(s,0)ds=-4\lambda_2\int_{-\Lambda_0/2}^{\Lambda_0/2}\int_{-1}^0\partial_s^2\psi_1(s,y)\psi_1(s,y)dsdy,
\end{equation*}
or
\begin{equation}\label{Okt31cc}
-\frac{\lambda_2}{\tau_*}(\sinh(2\tau_*)-2\tau_*)=\int_{-\Lambda_0/2}^{\Lambda_0/2}(2\theta F_1-F_2)\eta_1(s)ds.
\end{equation}

Let us evaluate $2\theta F_1-F_2$. First we note that
$$
2\theta\psi_{1y}(s,0)\eta_2-2\psi_{1yy}(s,0)\eta_2=2a\eta_2\tau_*^2\frac{\cos\tau_*s}{\sinh\tau_*},
$$
$$
\eta_1^2(\theta\psi_{1yy}(s,0)-\psi_{1yyy}(s,0))=0
$$
and
$$
2\eta_1(\theta\psi_{2y}(s,0)-\psi_{2yy}(s,0))=2\eta_1\Big(\theta\alpha_1+\alpha_2\cos(2\tau_* s)(\theta 2\tau_*\cosh(2\tau_*)-4\tau_*^2\sinh(2\tau_*))\Big).
$$
Further we have
\begin{equation*}
2\psi_{1s}(s,0)(\psi_{1sy}(s,0)\eta_1+\psi_{2s}(s,0))=2\psi_{1s}(s,0)\sin(2\tau_*s)(\frac{\tau_*^2}{2}\coth(\tau_*)-2\alpha_2\tau_*\sinh(2\tau_*))
\end{equation*}
and
\begin{eqnarray*}
&&2\psi_{1y}(s,0)(\psi_{2y}(s,0)+\psi_{1yy}(s,0)\eta_1)\\
&&=2\psi_{1y}(s,0)(\alpha_1+2\alpha_2\tau_*\cos(2\tau_*s)\cosh(2\tau_*)-\tau_*^2\cos^2\tau_*s).
\end{eqnarray*}
Therefore, we get
\begin{eqnarray*}
&&2\theta F_1-F_2=2a\eta_2\tau_*^2\frac{\cos\tau_*s}{\sinh\tau_*}+2\eta_1\Big(\theta\alpha_1+\alpha_2\cos(2\tau_* s)(\theta 2\tau_*\cosh(2\tau_*)-4\tau_*^2\sinh(2\tau_*))\Big)\\
&&-2\psi_{1s}(s,0)\sin(2\tau_*s)(\frac{\tau_*^2}{2}\coth(\tau_*)-2\alpha_2\tau_*\sinh(2\tau_*))\\
&&-2\psi_{1y}(s,0)(\alpha_1+2\alpha_2\tau_*\cos(2\tau_*s)\cosh(2\tau_*)-\tau_*^2\cos^2\tau_*s).
\end{eqnarray*}
This implies
\begin{eqnarray*}
&&\int_{-\Lambda_0/2}^{\Lambda_0/2}(2\theta F_1-F_2)\eta_1(s)ds=-\frac{\tau_*^2\Lambda_0}{\sinh^2\tau_*}\Big(\beta_1+\frac{\beta_2}{2}\Big)+\theta\alpha_1\Lambda_0\\
&&+\frac{\alpha_2\Lambda_0}{2}(\theta 2\tau_*\cosh(2\tau_*)-4\tau_*^2\sinh(2\tau_*))-\frac{\tau_*\Lambda_0}{2}(\frac{\tau_*^2}{2}\coth(\tau_*)-2\alpha_2\tau_*\sinh(2\tau_*))\\
&&+\tau_*\Lambda_0\coth\tau_*(\alpha_1+\alpha_2\tau_*\cosh(2\tau_*)-\frac{3}{4}\tau_*^2)
\end{eqnarray*}
or equivalently
\begin{eqnarray}\label{Nov3aa}
&&\frac{1}{\Lambda_0}\int_{-\Lambda_0/2}^{\Lambda_0/2}(2\theta F_1-F_2)\eta_1(s)ds=-\frac{\tau_*^2}{\sinh^2\tau_*}\Big(\beta_1+\frac{\beta_2}{2}\Big)+2\theta\alpha_1\nonumber\\
&&+\alpha_2(2\theta \tau_*\cosh(2\tau_*)-\tau_*^2\sinh(2\tau_*))-\theta\tau_*^2.
\end{eqnarray}

\subsection{Sign of $\lambda_2$}

Using (\ref{Nov3aa}), we write the equation (\ref{Okt31cc}) as
\begin{equation}\label{Nov3a}
-\frac{\lambda_2}{\tau_*\Lambda_0}(\sinh(2\tau_*)-2\tau_*)=f(\tau_*),
\end{equation}
where
\begin{equation}\label{Nov3b}
f(\tau)=-\frac{\tau_*^2}{\sinh^2\tau_*}\Big(\beta_1+\frac{\beta_2}{2}\Big)+2\theta\alpha_1
+\alpha_2(2\theta \tau_*\cosh(2\tau_*)-\tau_*^2\sinh(2\tau_*))-\theta\tau_*^2.
\end{equation}
The quantities $\theta$, $\alpha_j$ and $\beta_j$, $j=1,2$, are evaluated by the following formulas
$$
\theta=\nu(\tau_*),
$$
$$
2\alpha_1+2\theta\beta_1+\frac{1}{2}(\theta^2-\tau_*^2)=0,
$$
$$
4\tau_*\alpha_2\cosh 2\tau_*+2\theta\beta_2+\frac{1}{2}(\theta^2-3\tau_*^2)=0,
$$
$$
2(\theta-1)\beta_1+\theta+\frac{1}{2}(\theta^2-\tau_*^2)=0,
$$
and
$$
(2\theta-2\nu(2\tau_*))\beta_2+\theta\nu(2\tau_*)+\frac{1}{2}(\theta^2-3\tau_*^2)=0.
$$
Thus $f$ is well defined function of $\tau$ and it can be used for study of the sign of $\lambda_2$.

Evaluating the root $\tau_0$ of the equation $f(\tau)=0$ we get
\begin{equation}\label{J29ba}
\tau_0\approx 1.992.
\end{equation}
Therefore if $\tau\in (1,\tau_0)$ then
$$
\Lambda_2>0.
$$
According to (\ref{J29b}) and (\ref{J29ba}), we conclude that
\begin{equation}\label{J20c}
\mu_2>0\;\;\mbox{if $F<F_0$, where $F_0\approx 1.399 $}.
\end{equation}


\subsection{Upper  estimates of the Froude number and the validity of {\bf Assumption}}

  The following relation connected $d_-$and $d_+$ can be found in Sect. 2.1, \cite{KN11a} (see the formula (14) there):
\begin{equation}\label{Ju9a}
\frac{d_+}{d_-}=\frac{1+\sqrt{1+8d_-^3}}{4d_-^3},
\end{equation}
which implies
$$
d_+=\frac{1+\sqrt{1+8d_-^3}}{4d_-^2}.
$$
A necessary condition for existence of solitary wave is the lower estimate $F>1$, Therefore the depth $d=d_-$ corresponds to solitary waves and the corresponding Froude number is
$$
F=d_-^{-3/2}.
$$
Since $\theta=d_+^3$ we obtain the formula
\begin{equation}\label{J29bz}
\theta=\Big(\frac{F+\sqrt{F^2+8}}{4}\Big)^3F.
\end{equation}
Furthermore
$$
\tau\coth\tau=\theta.
$$
These two formulas give a connection between $\tau$ and $F$. In particular they imply the relation between $\tau_0$ and $F_0$ given by (\ref{J29ba}) and (\ref{J20c}).

The best known upper estimate for the Froude number is $F^2<2$. It can be derived from  \cite{Star} as it is explained in Introduction of \cite{We}. Therefore
 $F^2>2$  implies non-existence of solitary waves. Unfortunately $F_0<\sqrt{2}$ and it cannot be used for estimating of an interval for $F$, where the {\bf Assumption} is valid.
Another numerical estimate $F<1,29$ is obtained in \cite{LHF74}. This estimate together with (\ref{J20c}) leads to the relation  (\ref{J29da}).

\section{Acknowledgments}

I want to thank the anonymous referee for his excellent work as a referee of this paper. In particular, he has discovered  a mistake in numerical part of this paper, which was important for the final result.

\section{References}

{

\end{document}